\documentclass[12pt,a4paper]{amsart}
\usepackage[utf8]{inputenc}	
\usepackage[T1]{fontenc}
\usepackage[in,headings]{fullpage}
\usepackage{amsbsy, amscd, amsfonts, amsmath, amsrefs, amssymb, amsthm}
\usepackage{graphicx}
\usepackage{float}
\usepackage{color}   
\usepackage[colorlinks=true,linkcolor=blue,citecolor=blue,urlcolor=blue]{hyperref}

\newtheorem{theorem}{Theorem}

\newtheorem{proposition}[theorem]{Proposition}
\newtheorem{lemma}[theorem]{Lemma}

\theoremstyle{definition}

\theoremstyle{remark}

\newcommand{\R}{\mathbf{R}}

\renewcommand{\Re}{\mathop{\mathrm{Re}}\nolimits}

\newcommand{\Rzeta}{\mathop{\mathcal R }\nolimits}

\newfont{\cmbsy}{cmbsy10}
\newfont{\cmmib}{cmmib10}

\begin{document}

\title{Trivial zeros of Riemann auxiliary function.}
\author[Arias de Reyna]{J. Arias de Reyna}
\address{%
Universidad de Sevilla \\ 
Facultad de Matem\'aticas \\ 
c/Tarfia, sn \\ 
41012-Sevilla \\ 
Spain.} 

\subjclass[2020]{Primary 11M06; Secondary 30D99}

\keywords{función zeta, Riemann's auxiliary function}


\email{arias@us.es, ariasdereyna1947@gmail.com}


\begin{abstract}
It is proved that $s=-2n$ is a simple zero of $\Rzeta(s)$ for each integer $n\ge1$. Here $\mathop{\mathcal R}(s)$ is the function found by Siegel in Riemann's posthumous papers.

\end{abstract}

\maketitle

\section{Introduction}
The function 
\[\Rzeta(s)=\int_{0\swarrow1}\frac{x^{-s} e^{\pi i x^2}}{e^{\pi i x}-
e^{-\pi i x}}\,dx\]
was considered by Riemann and introduced by Siegel in \cite{Siegel}. The function $\Rzeta(s)$ is connected to the Riemann zeta function and is useful for the study of its zeros. 

In \cite{A166} we make a summary of the first properties of $\Rzeta(s)$, in particular it is shown that $\Rzeta(-2n)=0$ for each integer $n\ge1$. 
In \cite{A98}*{Thm.~7} it was shown that $s=-2n$ is a simple zero of $\Rzeta(s)$ for 
$n\ge 8495$ (i.~e.~$2n+1\ge 5408\pi$). Also, it is proved there that they are the only zeros on a region $G$ that contains most of the second and third quadrants.

Our objective here is to show that these \emph{trivial} zeros are all simple zeros.
The proof reduces to prove that  certain integrals related to the derivative $\Rzeta'(-2n)$ are not null. It is an exercise in estimation of integrals.

\section{An expression for the derivative.}

\begin{proposition}
At any trivial zero $s=-2n$, we have 
\begin{equation}\label{E:main}
\Rzeta'(-2n)=\frac{\omega\sqrt{\pi n}}{2}\Bigl(\frac{i n}{\pi e}\Bigr)^n\int_0^\infty\frac{(ex^2e^{-x^2})^n}{\sin(\omega x\sqrt{\pi n})}\,dx,
\end{equation}
where $\omega=e^{\pi i/4}$.
\end{proposition}
\begin{proof}
Start with \cite{A166}*{Prop.~6}, let $\omega=e^{\pi i/4}$ and $s=\sigma+it$ with $\sigma<0$ and $t\in\R$, then 
\begin{equation}\label{E:anterior}
\Rzeta(s)=\omega e^{\pi i s/4}\sin\frac{\pi
s}{2}\int_0^{+\infty}\frac{y^{-s}e^{-\pi y^2}}{\sin\pi\omega y}\,dy.
\end{equation}
From \eqref{E:anterior} it is clear that $\Rzeta(-2n)=0$ for an integer $n\ge1$, because $\sin\frac{\pi s}{2}$ vanish and the integral is convergent.

Differentiating \eqref{E:anterior} at $s=-2n$, we get 
\begin{equation}
\Rzeta'(-2n)=i^n\frac{\pi\omega}{2}\int_0^\infty\frac{y^{2n}e^{-\pi y^2}}{\sin\pi\omega y}\,dy.
\end{equation}
Changing the variable $\sqrt{\pi}\;y=\sqrt{n}\; x$ we get \eqref{E:main}.
\end{proof}

Our objective is to prove that the integrals in \eqref{E:main} do not vanish. The integral is a typical case for the Laplace approximation \cite{O}. In fact, the method
says that for $n\to\infty$ we have 
\begin{equation}\label{E:int}
\int_0^\infty\frac{(ex^2e^{-x^2})^n}{\sin(\omega x\sqrt{\pi n})}\,dx\sim \Bigl(\frac{\pi}{2n}\Bigr)^{1/2}\frac{1}{\sin(\omega \sqrt{\pi n})}.
\end{equation}
Hence, the integral does not vanish for $n\ge n_0$. But we want it to be $\ne0$ for all $n\ge1$. To get the asymptotic expansion with exact bounds is not so easy,and for small n may not be enough;  therefore, we will get only an inferior bound for its modulus.

\begin{theorem}
For all integers $n\ge1$ we have $\Rzeta'(-2n)\ne0$. Therefore, all trivial zeros are simple.
\end{theorem}
\begin{proof}
Given the approximate value of the integral in \eqref{E:int} and that $\sin(\omega \sqrt{\pi n})\approx \frac{i}{2} e^{-i\omega\sqrt{\pi n}}$, we will consider 
\begin{equation}\label{D:In}
I_n:= i e^{-i\omega\sqrt{\pi n}}\int_0^\infty\frac{(ex^2e^{-x^2})^n}{\sin(\omega x\sqrt{\pi n})}\,dx.
\end{equation}
Hence $I_n\sim 2\sqrt{\pi/2n}$ is almost a real number. By \eqref{E:main} and \eqref{D:In}
we have 
\[\Rzeta'(-2n)=\frac{\omega\sqrt{\pi n}}{2}\Bigl(\frac{i n}{\pi e}\Bigr)^n\frac{I_n}{ie^{-i\omega\sqrt{\pi n}}}.\]
To show that $\Rzeta'(-2n)\ne0$ we only have to show that $\Re I_n>0$. 

The function $ex^2e^{-x^2}$ has a maximum $=1$  at the point $x=1$. The main part of the integral is  for $x\in (1-n^{-1/2}, 1+n^{-1/2})$. Therefore, 
\[\Re I_n\ge A_n-B_n-C_n,\]
where 
\begin{align*}
A_n&=\Re\int_{1-n^{-1/2}}^{1+n^{-1/2}}\frac{i e^{-i\omega\sqrt{\pi n}}}{\sin(\omega x\sqrt{\pi n})}(ex^2e^{-x^2})^n\,dx\\
B_n&=\Bigl|\int_0^{1-n^{-1/2}}\frac{e^{-i\omega\sqrt{\pi n}}}{\sin(\omega x\sqrt{\pi n})}(ex^2e^{-x^2})^n\,dx\Bigr|\\
C_n&=\Bigl|\int_{1+n^{-1/2}}^\infty\frac{e^{-i\omega\sqrt{\pi n}}}{\sin(\omega x\sqrt{\pi n})}(ex^2e^{-x^2})^n\,dx\Bigr|.
\end{align*}
In Lemmas \ref{L:A}, \ref{Bn} and \ref{L:C} we show that for $n\ge9$ we have
\[A_n\ge \frac{1.471}{\sqrt{n}},\quad B_n\le \frac{1}{3\sqrt{n}},\quad
C_n\le \frac{1}{16\sqrt{n}}.\]
This proves our Theorem for all $n\ge9$. We end the proof by numerically computing the integral in \eqref{E:int} to verify that they are not null for $1\le n\le 9$. 
\end{proof}

\section{Bounds for \texorpdfstring{$A_n$}{An}, \texorpdfstring{$B_n$}{Bn}, and \texorpdfstring{$C_n$}{Cn}.}

\begin{lemma}\label{L:simple} We have
\begin{itemize}
\item[(a)] \label{I:a}
$\displaystyle{(ex^2e^{-x^2})^n\ge  \frac{2}{3} e^{-2n(x-1)^2}, \qquad 1-n^{-1/2}\le x\le 1+n^{-1/2}},\quad n\ge9$.
\item[(b)] \label{I:b} For $0<x<1$ we have $ex^2e^{-x^2}\le e^{-2(x-1)^2}$.
\item[(c)] \label{I:c} For $x\ge1$ we have $ex^2e^{-x^2}\le e^{-(x-1)^2}$.
\end{itemize}
\end{lemma}
\begin{proof}
The function $f(x):=2(x-1)^2+1-x^2+2\log x$ increases for $x>0$ since its derivative is $f'(x)=2x+2/x-4=2(x-1)^2/x\ge0$. Assertion (a) is equivalent to the inequality
$f(x)\ge \frac{\log(2/3)}{n}$ for $1-\frac{1}{\sqrt{n}}<x<1+\frac{1}{\sqrt{n}}$ and $n\ge9$. Hence, it is equivalent to $f(1-\frac{1}{\sqrt{n}})\ge \frac{\log(2/3)}{n}$. 
Expanding the logarithm, this is equivalent to 
\[2\sum_{k=3}^\infty \frac{1}{k n^{k/2}}\le \frac{\log(3/2)}{n}.\]
But 
\[2n\sum_{k=3}^\infty \frac{1}{k n^{k/2}}\le \frac{2n}{3}\sum_{k=3}^\infty \frac{1}{ n^{k/2}}= \frac{2n}{3}\frac{n^{-3/2}}{1-n^{-1/2}}\le \frac{2}{3}\frac{n^{-1/2}}{1-n^{-1/2}}\le \frac13\le\log(3/2).\]
This proves (a).

Taking logarithms  (b) is equivalent to $1+2\log x-x^2 +2(x-1)^2\le 0$.  Since $1+2\log x-x^2 +2(x-1)^2$ is increasing and takes the value $0$ at $x=1$, the result follows.

To prove (c), note that it is equivalent to $1+\log x\le x$ for $x\ge1$. But $x-1-\log x$ increases for $x>1$ and vanishes for $x=1$. 
\end{proof}

\begin{lemma}\label{L:A}
We have 
\begin{equation}
A_n\ge \frac{1.471}{\sqrt{n}},\qquad n\ge9.
\end{equation}
\end{lemma}
\begin{proof}
The factor in the integrand of $A_n$ is equal to 
\[-ie^{i\omega\sqrt{\pi n}}\sin(\omega x\sqrt{\pi n})=\frac{1}{2}
e^{i\omega\sqrt{\pi n}}(e^{-i\omega x\sqrt{\pi n}}-e^{i\omega x\sqrt{\pi n}})=
\frac{1}{2}(e^{-i\omega (x-1)\sqrt{\pi n}}-e^{i\omega (x+1)\sqrt{\pi n}})\]
Let $a_n=\sqrt{\pi n/2}$, then  $\omega \sqrt{\pi n}=(1+i)\sqrt{\pi n/2}=a_n(1+i)$.
\begin{align*}
-i e^{i\omega\sqrt{\pi n}}\sin(\omega x\sqrt{\pi n})&=
\frac{1}{2}(e^{-ia_n(1+i)(x-1)}-e^{ia_n(1+i)(x+1)})\\
&=\frac{1}{2}(e^{a_n(x-1)}e^{-ia_n(x-1)}-e^{-a_n(x+1)}e^{ia_n(x+1)})
\end{align*}
Consequently, for $1-n^{-1/2}\le x\le 1+n^{-1/2}$ and $n\ge9$
\[\Re\frac{i e^{-i\omega\sqrt{\pi n}}}{\sin(\omega x\sqrt{\pi n})}=\Re\frac{2}{e^{a_n(x-1)}e^{-ia_n(x-1)}-e^{-a_n(x+1)}e^{ia_n(x+1)}}=\Re
\frac{2e^{-a_n(x-1)}e^{ia_n(x-1)}}{1-e^{-2a_nx}e^{2ia_nx}}\]
Let $U=2e^{-a_n(x-1)}e^{ia_n(x-1)}$, $V=e^{-2a_nx}e^{2ia_nx}$. We have
$|V|= e^{-2a_nx}\le e^{-2a_n+2a_n n^{-1/2}}\le e^{-2\sqrt{2\pi}}$ for $n\ge9$,  and $|U|=2e^{-a_n(x-1)}\le 2e^{a_n n^{-1/2}}=2e^{\sqrt{\pi/2}}\le e^{\sqrt{2\pi}}$ therefore
\[\Re\frac{i e^{-i\omega\sqrt{\pi n}}}{\sin(\omega x\sqrt{\pi n})}=\Re(U+UV+UV^2+\cdots)\ge 2 e^{-a_n(x-1)}\cos(a_n(x-1))-e^{\sqrt{2\pi}}\frac{e^{-2\sqrt{2\pi}}}{1-e^{-2\sqrt{2\pi}}}.\]
To simplify notation, put
\[f(x)=2 e^{-a_n(x-1)}\cos(a_n(x-1))-\frac{e^{-\sqrt{2\pi}}}{1-e^{-2\sqrt{2\pi}}}.\]
Then
\[A_n\ge\int_{1-n^{-1/2}}^{1+n^{-1/2}}f(x)(ex^2e^{-x^2})^n\,dx,\]
and by Lemma \ref{L:simple}, {(a)} 
\[A_n\ge\frac{2}{3}\int_{1-n^{-1/2}}^{1+n^{-1/2}}f(x)e^{-2n(x-1)^2}\,dx\]
Changing variables $\sqrt{2n}(x-1)=y$, we get
\[A_n\ge\frac{2}{3\sqrt{2n}}\int_{-\sqrt{2}}^{\sqrt{2}}\Bigl(2 e^{-\tfrac12y\sqrt{\pi}}\cos(\tfrac12y\sqrt{\pi})-\frac{e^{-\sqrt{2\pi}}}{1-e^{-2\sqrt{2\pi}}}\Bigr)e^{-y^2}\,dy\]
Computing numerically this integral, we get 
\[A_n\ge \frac{1.471}{\sqrt{n}}.\qedhere\]
\end{proof}
\begin{lemma}\label{L:B}
We have 
\begin{equation}\label{Bn}
B_n\le \frac{1}{3\sqrt{n}},\qquad n\ge9.
\end{equation}
\end{lemma}
\begin{proof}
Since we integrate for $x>0$ we have
\[|e^{-i\omega\sqrt{\pi n}}|=e^{\sqrt{\pi n/2}},\quad |\sin(\omega x\sqrt{\pi n})|
\ge \tfrac12(e^{x\sqrt{\pi n/2}}-e^{-x\sqrt{\pi n/2}})=\sinh(a_n x).\]
It follows that
\begin{equation}\label{E:partial}
B_n\le  e^{a_n}\int_0^{1-n^{-1/2}}\frac{(ex^2e^{-x^2})^n}{\sinh(a_n x)}\,dx,\quad 
C_n\le  e^{a_n}\int_{1+n^{-1/2}}^\infty\frac{(ex^2e^{-x^2})^n}{\sinh(a_n x)}\,dx.
\end{equation}
The function $(ex^2e^{-x^2})$ increases for $0<x<1$ and is less than $\frac13$ for $x=\frac13$. Therefore,
\[B_n\le \frac{e^{a_n}}{3^{n-1}}\int_0^{1/3}\frac{ex^2e^{-x^2}}{\sinh(a_n x)}\,dx+
e^{a_n}\int_{1/3}^{1-n^{-1/2}}\frac{(ex^2e^{-x^2})^n}{\sinh(a_n x)}\,dx.\]
For the first integral, we have (for $n\ge9$)
\[=\frac{e^{a_n}}{3^{n-1}a_n}\int_0^{1/3}\frac{a_n x}{\sinh(a_n x)}exe^{-x^2}\,dx\le \frac{e^{a_n}}{3^{n-1}a_n}\int_0^{1/3}exe^{-x^2}\,dx< \frac{e^{a_n}}{3^{n}a_n}\le \frac{1}{574\sqrt{n}}\]

By Lemma \ref{L:simple} {(b)}, we have 
\[e^{a_n}\int_{1/3}^{1-n^{-1/2}}\frac{(ex^2e^{-x^2})^n}{\sinh(a_n x)}\,dx\le e^{a_n}\int_{1/3}^{1-n^{-1/2}}\frac{e^{-2n(x-1)^2}}{\sinh(a_n x)}\,dx\]
The function $e^x/\sinh x$ decreases for $x>0$ with limit $2$ for $x\to+\infty$
Therefore, for $n\ge9$ and $x>1/3$ we have 
\begin{equation}\label{sinh}
\frac{1}{\sinh(a_n x)}\le \frac{e^{a_9 /3}}{\sinh(a_9/3)}e^{-a_n x}\le \frac{11}{5}e^{-a_n x}.
\end{equation}
Hence, 
\[e^{a_n}\int_{1/3}^{1-n^{-1/2}}\frac{(ex^2e^{-x^2})^n}{\sinh(a_n x)}\,dx\le \frac{11e^{a_n}}{5}\int_{1/3}^{1-n^{-1/2}}e^{-2n(x-1)^2-a_nx}\,dx\]
\[\le\frac{11e^{a_n}}{5\sqrt{2n}}\int_{-\infty}^{-\sqrt{2}}e^{-y^2-y\sqrt{\pi}/2-a_n}\,dy < \frac{1}{3\sqrt{n}}.\]
Not only this integral, the bound $\frac{1}{3\sqrt{n}}$ is also true when we include the  term $\frac{1}{574\sqrt{n}}$. So, we have \eqref{Bn}.
\end{proof}

\begin{lemma}\label{L:C}
We have 
\begin{equation}\label{Cn}
C_n\le \frac{1}{16\sqrt{n}},\qquad n\ge9.
\end{equation}
\end{lemma}
\begin{proof}
By \eqref{E:partial} and then by \eqref{sinh}
\[C_n\le  e^{a_n}\int_{1+n^{-1/2}}^\infty\frac{(ex^2e^{-x^2})^n}{\sinh(a_n x)}\,dx
\le \frac{11e^{a_n}}{5}\int_{1+n^{-1/2}}^\infty(ex^2e^{-x^2})^n e^{-a_n x}\,dx.\]
Applying Lemma \ref{L:simple} {(c)} we get 
\[C_n\le \frac{11e^{a_n}}{5}\int_{1+n^{-1/2}}^\infty e^{-n(x-1)^2-a_n x}\,dx=
\frac{11e^{a_n}}{5\sqrt{n}}\int_{1}^\infty e^{-y^2-y\sqrt{\pi/2}-a_n}\,dy\le \frac{1}{16\sqrt{n}}.\qedhere\]
\end{proof}

\end{document}